\documentclass{amsart}
\usepackage[T1]{fontenc}
\usepackage[utf8]{inputenc}
\usepackage{amsaddr}
\usepackage[T1]{fontenc}
\usepackage[margin=1.1in]{geometry}
\usepackage{graphicx}
\usepackage{hyperref}
\usepackage{color}

\usepackage{amssymb,amsfonts,amsmath,latexsym}
\usepackage{booktabs}
\usepackage{multirow}
\usepackage{cellspace} %
\usepackage{wrapfig}
\usepackage{subfigure}
\usepackage{caption}
\usepackage{xcolor}
\usepackage{enumitem}
\setlist[enumerate]{label*=\arabic*.}
\usepackage[noabbrev]{cleveref}  
\usepackage{algorithmic}
\usepackage{longtable}
\usepackage{array} 
\usepackage[ruled,lined,linesnumbered]{algorithm2e}

\usepackage[
    backend=biber,
    style=numeric,
    giveninits=true,
    maxbibnames=999,
    maxcitenames=999,
    sortcites=true
]{biblatex}

\newcommand{\bm}[1]{\mathbf{#1}}
\newcommand{\ttt}[1]{\texttt{#1}}

\newcommand{\R}{{\mathbb{R}}}

\newcommand{\cI}{{\mathcal{I}}}
\newcommand{\cS}{{\mathcal{S}}}

\newcommand{\bell}{\boldsymbol{\ell}}
\newcommand{\pstredtools}{\texttt{\small ps3dtools}}
\DeclareMathOperator{\diag}{diag}

\title[A Matlab/Octave package for photometric stereo 3D reconstruction]{A 
Matlab/Octave package for photometric stereo 3D reconstruction}

\author[A. Buccini, E. Crabu, F. Pes, G. Rodriguez]{
Alessandro Buccini\textsuperscript{1}
\and Elisa Crabu\textsuperscript{1}
\and Federica Pes\textsuperscript{2}
\and Giuseppe Rodriguez\textsuperscript{1}
}

\keywords{shape from shading, photometric stereo, 3D reconstruction, non 
Lambertian surfaces}

\subjclass[2020]{68U10, 65D19, 65Y15}

\begin{document}
\maketitle

\begin{center}
\textsuperscript{1} Department of Mathematics and Computer Science,
University of Cagliari, via Ospedale 72, 09124 Cagliari, Italy\\
\texttt{\{alessandro.buccini, elisa.crabu, rodriguez\}@unica.it}

\medskip

\textsuperscript{2} Department of Science and High Technology,
University of Insubria, via Valleggio 11, 22100 Como, Italy\\
\texttt{federica.pes@uninsubria.it}
\end{center}

\begin{abstract}
The photometric stereo model allows the construction of 3D virtual replicas of 
real objects from images acquired under different light conditions.
The assumptions of this model are strongly restrictive and are seldom 
fully satisfied in applications.
Indeed, the surface of the observed objects
may deviate from ideal conditions, and 
the geometry of the scene is often only an approximation of the one assumed by
the model. 
Additionally, real-world surfaces may present reflection phenomena, which 
introduce distortions in reconstructions.
In this paper, we present 
a revised version of the  Matlab package \ttt{\small ps3d}~\cite{cdfrv21},
developed to implement the photometric stereo model under various
conditions. 
The new package contains various extensions, the most notable one is the
presence of an alternative model to address the presence of a reflective object
surface.
In addition, it implements particular metrics to analyze the compliance of a
dataset to model assumptions, and it is compatible with GNU Octave.
The tool allows users to test the model under different
observation configurations, also in presence of reflection, using
either synthetic or experimental data.
\end{abstract}

\section{Introduction}

The problem of reconstructing a three-dimensional replica of a given object,
starting from easily obtainable data, has been studied over years through
several different techniques. 
For example, \emph{normal integration}~\cite{brumentrnb2024, cao2022,
cao2024, qda18} estimates the shape of an object from its normal map,
while \emph{reflectance transformation imaging
(RTI)}~\cite{palma2010,diaz2015,lewis2024} constructs a 2D interactive model
that derives all possible lighting of the target on the basis of 
a large amount of images with different light sources. 
Several deep-learning methods have also been
developed for this purpose; see, e.g.,~\cite{huang2024, xiong2022, li2022}.

\emph{Shape from shading (SfS)}~\cite{horn1989, klette1998, slabaugh2001,
dfs2008, rpjm1999, YaBr2020, tf2016} techniques address the problem of
creating a 3D digital model of a surface from a collection of images. 
Among the best known approaches for SfS are 
multiview~\cite{dyer2001, jin2008, langguth2016} and 
photometric stereo (PS)~\cite{woodham1979, woodham1980}.
In the first one, the lighting of the scene is fixed and the observer, i.e., the
camera, moves through it collecting a sequence of pictures.
On the contrary,
PS considers a fixed camera pointing at an object of interest, and each image
correspond to a different light source.
In this paper we focus on this approach.

A PS dataset is composed of images collected under
different lighting conditions by a fixed camera. The standard model, also known 
as
\emph{Lambert's law}, assumes that
the observed light intensity is independent of the point of view, and depends 
exclusively on the relative position between the observed object and the light 
source.
Such property of a surface is called \emph{Lambertianity}.
Another important assumption is that the light source is at
infinite distance from the observed scene,  so that all light rays are parallel
to a known direction.
For simplicity, we consider the camera positioned at a large distance from the
scene, which corresponds to orthographic projection.
If needed, perspective parameters can be determined by camera calibration
techniques.

The PS problem has been handled in various settings~\cite{khanian2018,
Radow2018, barsky2003}. In~\cite{kozera1991}, Kozera proposed a formulation
that leads to a Hamilton-Jacobi problem, later numerically analyzed
in~\cite{mecca2011, mecca2013}, for which at least two images are needed. 
An alternative formulation, described for example in~\cite{cdfrv21}, 
requires at least three images with different lighting, and
determines the solution of the PS problem in two steps.
Initially, the normal vector field is estimated from data, then these vectors
are used to solve a Poisson problem to determine an approximation of the shape
of the object. 
PS has been employed in many research areas, from medical
studies~\cite{sun2008, caops2022} to archaeology~\cite{brument2024,melou2022,
ikeuchi2022}. 

The basic assumption of known light positioning is very
limiting, since in real-world applications
it is not always possible to measure the exact light location in the scene
during the acquisition of the images. 
Hayakawa~\cite{hayakawa1994} showed that, if at least six differently lighted
images are available, then the lights position can be estimated from the
dataset. This procedure has been used in~\cite{cdfrv21} to create a practical 
algorithm which determines lighting directions from data, solving at the same
time the so-called \emph{bas-relief} ambiguity~\cite{basrelief99}.

Although the method produces accurate results when applied to synthetic
datasets, numerical tests show that the computation with experimental data
often suffers of some issues. 
In some cases, the model returns inaccurate reconstruction, or the algorithm
may even  break down before determining a solution.
These issues arise from ``non-ideal'' dataset, due to unmet model assumptions.
To handle such cases, in~\cite{cprt23} two procedures for dataset
pre-processing have been proposed, based on numerical indicators of
ideality. These methods lead to a reduction of a dataset that may
improve the solution quality, allowing to determine
reasonable reconstructions. 

Another important issue of the standard PS model is that real objects rarely 
behave as Lambertian surfaces, and reflect (part of) the light in the
direction of the observer, making the observed intensity depend on the point
of view. This results in distorted data information.
The non-Lambertian case has been extensively studied: some learning-based
method are described in~\cite{logothetis2021, ju2022}, while
in~\cite{ikehata2012} a statistical approach, based on solving a regression
problem, is used to estimate the normal vector field. Differently,
in~\cite{wu2011}, starting from the observation that non-Lambertian surfaces
are subject of shadows and reflection, a solution method based on low-rank
matrix factorization affected by sparse errors has been proposed.
A pre-processing technique to treat the same condition has been studied
in~\cite{rrmyb23}.

In this paper, we present a new method to treat the reflection issue.
The problem is re-formulated as a minimization problem and solved using
the Major\-ization-Minimization algorithm~\cite{huang2017}. 
The model is defined under the assumption of known lighting, 
and here we present preliminary results
obtained from the application to synthetic datasets. 

Since an implementation of PS related algorithms is not trivial,
we present in this work a revised version of the Matlab package
\ttt{\small  ps3d} from~\cite{cdfrv21}.
It extends the original software by including all the
configurations described above and improving some routines.
The user has the possibility to test the algorithms using both
synthetic and experimental datasets, through easily modifiable
demonstration scripts. 
The new package has been renamed to \pstredtools\ to
avoid confusion with the name of one of its functions.

Several implementations of the photometric stereo model are publicly available.
Some are based on classical Lambertian approaches, such as the Python
implementation of Woodham's algorithm provided in~\cite{visiont3lab} and the
Matlab implementation~\cite{niklassanden}, which follows the formulation of
Hayakawa. More recently, LINO-UniPS~\cite{LINO} tackles the photometric stereo
setting, where lighting conditions are arbitrary and unknown, by means of
learning-based methods designed to estimate surface normals without explicit
knowledge of light directions. However, none of these tools provides an approach
to assess the conformity of a dataset with the model assumptions, nor do they
include an explicit model for reflective surfaces. The package \pstredtools\ 
presented in this paper addresses all these aspects.

The paper is organized as follows. 
\Cref{sec:ps} describes the PS model in
detail, while non-ideality in data is treated in \Cref{sec:ideal}, and
reflection phenomena in \Cref{sec:reflmod}. \Cref{sec:package} illustrates the 
\pstredtools\ package, showing some numerical experiments that can be
performed with it.
Finally, \Cref{sec:con} summarizes considerations related to the paper content.

\section{The photometric stereo model}\label{sec:ps}

In this section, we briefly recap the main features of the PS mathematical
model.
An extensive description can be found in~\cite{cdfrv21}.

The basic configuration considers an object placed at the origin of a reference
system in $\R^3$, and a camera positioned so that its optical axis coincides
with the $z$-axis and the lens is directed towards the target. We assume
the point of view at infinite distance, allowing for the use of
orthographic projection. The images used by the method are obtained moving a
light source, theoretically at infinite distance from the object, to different
positions in the upper half-space. For each position, a picture is collected.
 
Images are assumed to have resolution $(r +2) \times (s+2)$.
Each of them has sides of length $A$ and $B$, where $h =A/(r+1)$ is the size of
a pixel. Assuming square pixels, we set $B=(s+1)h$.
Therefore, each picture defines a domain $\Omega=[-A/2,A/2]\times[-B/2,B/2]$.
Assuming that the surface of the object is represented by the explicit function
$z = u(x,y)$, $(x,y)\in \Omega$, the normalized normal vector is given by 
$$
\bm{n} (x,y) = \frac{\begin{bmatrix}-u_x(x,y),& -u_y(x,y),& 
1\end{bmatrix}^T}{\sqrt{1+\| \nabla u(x,y) \|^2} },
$$
where $u_x$, $u_y$ denote the partial derivatives of $u$ with respect to $x$
and $y$, respectively, $\nabla u$ is the gradient of $u$, and $\| \cdot \|$ is
the Euclidean norm. 

We consider a uniform discretization the $\Omega$ domain, i.e., a grid of
points of coordinates $(x_i,y_j)$, with $x_i=-\frac{A}{2}+ih$, $i=0, \dots,
r+1$, and $y_j=-\frac{B}{2}+jh$, $j=0, \dots, s+1$.
As it is customary, we order the pixels lexicographically, 
so that $u_k=u(x_i,y_j)$, where
$k := (i-1)s + j = 1, \dots , p$ and $p = rs$ is the total number of 
grid points.

Let $q$ represent the number of light directions, and denote 
by $\bell_t =\left[\ell_{1,t},\; \ell_{2,t},\; \ell_{3,t}\right]^T$, $t=1, 
\dots, q$,
the vector that points from the object to the light source.
The light intensity is proportional to its norm 
that, for simplicity, is assumed
to be unitary, meaning that all the lights have the same intensity.

We assume that the surface is Lambertian, i.e., that it satisfies
\emph{Lambert's cosine law}
\begin{equation}\label{eqlam}
\rho(x,y) \langle \bm{n}(x,y), \bell_t \rangle = \cI_t(x,y), \quad t=1, \dots, 
q,
\end{equation}
where $\langle\cdot,\cdot\rangle$ denotes the inner product in $\R^3$.
The function $\rho(x,y)$ represents the \emph{albedo}, which takes into account
the partial absorption of the light from the surface, and $\cI_t(x,y)$ is the
pixel value at the point $(x,y)$ of the $t$th image.
Discretizing \cref{eqlam}, Lambert's law takes the matrix form
\begin{equation}\label{eqmatrix}
D N^T L = M,
\end{equation}
where 
$$
D=\diag(\rho_1,\ldots,\rho_p), \quad
N=[\bm{n}_1,\ldots,\bm{n}_p], \quad
L=[\bell_1,\ldots,\bell_q], 
$$
being $\rho_k$ and $\bm{n}_k$ the value of the albedo and the normal vector at
the $k$th point of the discretization, respectively. 
The matrix $L$ contains the light directions and
each columns of the $p\times q$ matrix $M$ is a vectorized image of the
dataset.

If $L$ is known, matrices $D$ and $N$ can be easily estimated 
from~\eqref{eqmatrix} by solving a least-squares problem or, equivalently, by
computing the pseudoinverse matrix $L^\dagger$.
It is then immediate to observe that this produces a unique solution only if
$q\geq 3$.

Under unknown lighting, that is, when the position of the light sources is
not available, the normal vector field can be estimated
from~\eqref{eqmatrix} using Hayakawa's procedure~\cite{hayakawa1994}, that we
now summarize.
Let us consider the compact singular value decomposition (SVD) of the data
matrix $M = U \Sigma V^T$, where $\Sigma =\diag(\sigma_1, \dots, \sigma_q)$ 
contains the singular values of $M$, $U \in \R^{p \times q}$ and 
$V \in \R^{q \times q}$ have orthonormal columns. 
If all the assumptions of the
model are met, the rank of $M$ is $3$, i.e.,
$\sigma_4=\cdots=\sigma_q=0$. 
However, due to experimental errors and non-ideality of data, in general, 
$\sigma_i$ does not vanish for $i>3$.

To restore the theoretical rank, we truncate the SVD letting
$$
W=[\sigma_1 \bm{u}_1, \sigma_2 \bm{u}_2, \sigma_3 \bm{u}_3], \qquad
Z= [\bm{v}_1, \bm{v}_2, \bm{v}_3]^T,
$$
and substituting the matrix $M=WZ$ to the original dataset.
Then, Hayakawa's procedure consists of determining the $3\times 3$ symmetric
positive definite matrix $G$ such that 
\begin{equation}\label{hayaprob}
\diag(Z^T GZ) = \bm{1},
\end{equation}
with $\bm{1} = (1, \dots, 1)^T \in \R^q$, and computing the Cholesky
factorization $G = R^TR$.
It is possible to see that from $R$ one can determine an approximation of both
the matrices $N$ and $L$, up to a unitary transformation; see~\cite{cdfrv21}
for details.
The solution of~\eqref{hayaprob} is determined by solving the a linear
least-squares problem with a coefficient matrix of size $q \times 6$. 
This shows that the solution is unique if $q \geq 6$, so that at least 6 images
are needed.
 
In both the known and unknown lighting scenarios, once $N$ is computed
the 3D profile of the object can be recovered, e.g., by
numerically differentiating the normal vector field and
solving the Poisson partial differential equation 
$$\Delta u(x,y) = f(x,y),$$
where $\Delta$ denotes the Laplace operator.
Suitable boundary conditions must be imposed to guarantee the unicity of the
solution; see~\cite{cdfrv21} for a discussion of this issue.

\section{Data ideality}\label{sec:ideal}

When processing experimental data, that not necessarily verify the ideality
assumptions of PS, the computational procedure may fail. Indeed,
light sources positioned too close to the target, or the presence of non
Lambertian surfaces, may introduce distortion in the data and result in
strongly inaccurate reconstructions.
Moreover, when Hayakawa's method is implemented to estimate the light
positions, the computation may even breakdown because of the above issues. 
The breakdown is caused by the fact that,
when the matrix $M$ contains non-ideal data, the symmetric matrix $G$, solution
of~\eqref{hayaprob}, may not be numerically positive definite, preventing the
computation of its Cholesky factor.

This kind of problems may happen in many experimental scenarios when data
is collected on field.
In particular, the requirement of positioning light sources at a large
distance from the target may be very limiting in real-world applications.

This was observed, e.g., in~\cite{cprt23}, where the aim was to reconstruct
archaeological engravings found in \emph{Domus de Janas} sites.
Such prehistoric remains are funerary hypogea formed by very small rooms
carved into stone, with partially decorated walls.
Due to the structure of the sites, the light sources have to be positioned
quite close to the target, failing to meet one of the model assumptions. 
Moreover, rock surfaces deviate from Lambertianity, due to their composition,
so that the basic mathematical model is only an approximation of reality.

In~\cite{cprt23}, it was observed that processing experimental datasets,
typically composed by 10--20 images with different lighting, lead to distorted
3D reconstructions.
At the same time, removing specific images from the dataset may have the effect
of dramatically improving the accuracy of the computed solution.
Such ``non-ideal'' images to be removed were initially chosen by chance, as
there is no visual characteristic that identifies them.

To help practitioners to reduce the number of images and produce a more
effective dataset, two numerical indicators for ``data ideality'' were
developed in~\cite{cprt23}.
Such procedures iteratively remove an image from a dataset, and check the
remaining data for consistency.
This results in a list of images to be excluded, that the user may use to
improve results.

The first indicator checks for the positivity of the matrix $G$ 
in~\eqref{hayaprob}, and selects the configuration that corresponds to the 
largest positive value of the minimal eigenvalue of $G$.

The second one starts from an alternative formulation of Hayakawa's procedure,
proposed for the first time in~\cite{cprt23}.
It can be shown that solving~\eqref{hayaprob} is equivalent to finding an upper
triangular matrix $R$ such that
\begin{equation}\label{nonlineq}
||R\bm{z}_t\|=1, \quad t=1,\ldots,q.
\end{equation}
This produces a nonlinear least-squares problem, which is solved by a
Gauss--Newton algorithm proposed in~\cite{pr20}.
It turns out that in the presence of distorted data, the Jacobian of the
problem becomes rank-deficient at convergence.
The second ideality indicator selects the configuration which ensures that 
the last singular value of the Jacobian is nonzero.

The two indicators have been implemented in \pstredtools\
toolbox. A wide experimentation is available in~\cite{cprt23}, and some results
are reported in \Cref{sec:package}.

\section{Dealing with reflective surfaces}\label{sec:reflmod}

As already remarked, the standard PS model does not perform well when
the observed surface behaves as a non-ideal Lambertian reflector, that is, when
it reflects light in the direction of the point of view.
As a result, some pixel may exhibit specular reflection, causing distortion in
reconstructions. 
To handle this issue, it is important to explicitly account for these distorted
pixels in the mathematical model. 

To address this problem, we modified the basic PS model.
This reformulation stands on the assumption of known lighting.
A generalization to unknown lighting will be treated in future work.

Let us consider the following modification of problem~\eqref{eqmatrix} 
\begin{equation}\label{newmtxlambert}
D N^T L = M+E,
\end{equation}
where $E$ collects reflection points. Since just few points typically present
reflection, we can assume that $E$ is a sparse matrix, i.e., that most of its 
entries vanish.
The matrix $L$ is assumed to be available.

We rewrite equation~\eqref{newmtxlambert} as
$$
L^T ND = (M + E)^T,
$$
or
\begin{equation}\label{eqprob}
L^T \widetilde{N} = \widetilde{M},
\end{equation}
with $\widetilde{M}=(M+ E)^T$ and $\widetilde{N}=ND$. 
Our objective is, given the measured data matrix $\widetilde{M}$ and the
lighting positions $L$, to determine $\widetilde{N}$. 

Let $0<\tau\leq 1$, we define the $\tau$-norm of a matrix $X  \in \R^{m \times 
n}$ as 
$$
\|X\|_\tau^\tau = \sum_{i=1}^m \sum_{j=1}^n \lvert x_{ij}\rvert^\tau.
$$
With a slight abuse of notation,
we refer to $\| \cdot \|_{\tau}$ as the $\ell^{\tau}$-norm of $X$ for all
$\tau$, even though it is not a norm when $\tau< 1$,
since it does not satisfy the triangular inequality.
Then, an approximate solution of~\eqref{eqprob} can be obtained by solving the
minimization problem 
\begin{equation}\label{reflprobl}
\begin{aligned}
\widetilde{N}^{*} & = \arg\min_{\widetilde{N}} \|\widetilde{M} - L^T 
\widetilde{N} \|_\tau^\tau\\
& = \arg\min_{\widetilde{N}} \sum_{i=1}^q \| \widetilde{\bm{m}}_i - L^T 
\widetilde{\bm{n}}_i\|_\tau^\tau, \quad 0<\tau\leq1,
\end{aligned}
\end{equation}
where $\widetilde{\bm{m}}_i$ and $\widetilde{\bm{n}}_i$ denote the $i$-th
columns of matrices $\widetilde{M}$ and $\widetilde{N}$, respectively.
In~\eqref{reflprobl} we consider $0<\tau \leq 1$ to enforce sparsity on the
residual $\widetilde{M} - L^T \widetilde{N}^{*}$ that approximates $E$; for 
details and
Bayesian justification see~\cite{buccini2020}. 

The $\ell^{\tau}$-norm is non-smooth and, when $\tau<1$, it is non-convex.
This is usually the case of interest. To remove the non-differentiability of
the $\ell^{\tau}$-norm we approximate it by the smooth function  
\begin{equation*}
\| \bm{x} \|_{\tau, \varepsilon}^\tau
:= \sum_{j=1}^n (x_j^2 + \varepsilon^2 )^{\tau/2} 
\approx \| \bm{x} \|_\tau^\tau ,
\qquad \bm{x} \in \R^n, 
\end{equation*}
where $\varepsilon >0$ is a user-defined parameter. Note that the introduction
of $\varepsilon$ ensures that $\| \bm{x} \|_{\tau, \varepsilon}^\tau$ is
everywhere differentiable, even when one component of $\bm{x}$ vanishes.
To compute an approximate solution of~\eqref{reflprobl}, we use the
Majorization-Minimization algorithm~\cite{huang2017}.

\subsection{The Majorization-Minimization method} 

We briefly sketch the Majorization-Minimization (MM) algorithm when applied to
a problem of the form
\begin{equation*}\label{illp}
A\bm{x} = \bm{b},
\end{equation*}
with $\bm{b}\in \R^m$, $\bm{x} \in \R^n$, and  $A \in \R^{m \times n}$. The 
extension to our case is trivial.

We introduce the $\ell^\tau$ minimization problem
\begin{equation}\label{minlplq}
\bm{x}^* = \arg\min_\bm{x} \frac{1}{\tau} \left\|A\bm{x} - \bm{b} 
\right\|_\tau^\tau,
\end{equation}
where $0 < \tau \leq 1$.

Let us observe that this model was developed to tackle
ill-posed inverse problem and is a special case of the so-called
$\ell^p$-$\ell^q$ regularization; see, e.g.,~\cite{BR23,huang2017,LMSR15}. A
Matlab implementation can be found in~\cite{BR24}. 

An approximate solution of~\eqref{minlplq} can be found by minimizing the
differentiable functional 
\begin{equation}\label{functional}
J_\varepsilon (\bm{x}) = \frac{1}{\tau} \|A\bm{x} - \bm{b} \|_{\tau, 
\varepsilon}^\tau, \quad \varepsilon>0,
\end{equation}
whose minimum $\bm{x}_\varepsilon^*$ is usually fairly close to $\bm{x}^*$. 
Problem~\eqref{functional} can be solved using the MM algorithm, an iterative
method that generates a sequence  $\{\bm{x}^{(k)}\}$ converging to a
stationary point of $J_\varepsilon$. 

Given an approximation $\bm{x}^{(k)}$ of $\bm{x}^*$, the
method first constructs the so-called \emph{adaptive quadratic tangent
majorant} $Q(\bm{x},\bm{x}^{(k)})$, a quadratic function that majorizes
$J_\varepsilon$ everywhere and is tangent to it in $\bm{x}^{(k)}$. Then, the
new iterate  $\bm{x}^{(k+1)}$ is obtained by minimizing
$Q(\bm{x},\bm{x}^{(k)})$.
More in detail, we wish to construct $Q(\bm{x},\bm{x}^{(k)})$ such that:
\begin{itemize}
\item
$Q\left(\bm{x}^{(k)},\bm{x}^{(k)}\right)=J_\varepsilon\left(\bm{x}^{(k)}\right),
 \quad \nabla Q\left(\bm{x}^{(k)},\bm{x}^{(k)}\right)=\nabla 
J_\varepsilon\left(\bm{x}^{(k)}\right)$, 
\item $Q\left(\bm{x},\bm{x}^{(k)}\right)\geq J_\varepsilon(\bm{x})$, 
$\forall \bm{x} \in \R^n$.
\end{itemize}

Let
\begin{equation*}
\bm{u}^{(k)}= A \bm{x}^{(k)} - \bm{b},
\end{equation*}
and define the weight matrix
\begin{equation}\label{Wk}
W^{(k)}= \diag \left( \left((\bm{u}^{(k)})^2 + 
\varepsilon^2\right)^{\frac{\tau}{2}-1}\right),
\end{equation}
where all the operations are meant element-wise.
The \emph{adaptive quadratic tangent majorant} $Q$ is given by
\begin{equation}\label{adfun}
Q(\bm{x}, \bm{x}^{(k)}) = \frac{1}{2} \left\|\left(W^{(k)}\right)^{\frac{1}{2}} 
(A\bm{x} -\bm{b})\right\|_2^2.
\end{equation}
Then, the new approximate solution is obtained by 
\begin{equation*}
\bm{x}^{(k+1)}= \arg \min_\bm{x} Q(\bm{x}, \bm{x}^{(k)}).
\end{equation*}
In~\cite{huang2017}, a projection in a generalized Krylov subspace was used to
reduce the computational cost when $A$ is very large. In our application this
is usually not necessary, therefore, we do not discuss it here.

We employ the MM algorithm above to solve~\eqref{reflprobl}.
At each iteration $k$ of the MM method, for each column of $N$ and $M$, we 
construct the vector  $\bm{u}_i^{(k)}=L^T \widetilde{\bm{n}}^{(k)}_i - 
\widetilde{\bm{m}}_i$. Then, the adaptive weight $W^{(k)}_i$ is computed as
\begin{equation*}
W^{(k)}_i = \diag \left( \left( \left(\bm{u}_i^{(k)}\right)^2 + 
\varepsilon^2\right)^{\frac{\tau}{2}-1}\right).
\end{equation*}
Therefore, the new iterate is
\begin{equation*}
\widetilde{\bm{n}}^{(k+1)}_i = \arg\min_{\widetilde{\bm{n}}}\left \| 
\left(W_i^{(k)}\right)^{\frac{1}{2}} \left((L^T \widetilde{\bm{n}} - 
\widetilde{\bm{m}}_i\right) \right\|_2^2.
\end{equation*}
After a sufficiently large number of iterations are executed for each column
of $\widetilde{N}$, an approximation of the normal vector field is found.

We outline the MM method in Algorithm~\ref{mmalgo}. We stop the process when
either a fixed number of iterations is performed, or when
\begin{equation}\label{stop}
\sum_{i=1}^q\left\| \bm{n}^{(k+1)}_i - \bm{n}^{(k)}_i \right\|_{2}^2 \leq 
\gamma^2 \sum_{i=1}^q\left\| \bm{n}_i^{(k)}\right\|_{2}^2
\end{equation}
for a given $\gamma>0$.
Observe that the internal \ttt{for} loop at Line~\ref{for} is suitable for
parallel execution. Moreover, the minimization in Line~\ref{mini} can be
performed using a QR factorization of $\left(W_i^{(k)}\right)^{\frac{1}{2}} 
L^T$ since the
involved matrices are not large.

\begin{algorithm}
\caption{The Majorization-Minimization algorithm.}
\label{mmalgo}

\KwIn{$\widetilde{M}$, $L^T$, \emph{maxit} (maximum number of iterations), 
$\varepsilon$, $\gamma$}

Set $\widetilde{N}$ a $3 \times p$ zero matrix\;

\For{$k = 1, \dots, \textit{maxit}$ \label{for}}{
    $\widetilde{N}_{\rm old} = \widetilde{N}$\;
    
    \For{$i = 1, \dots, q$}{
        $\bm{u}_i^{(k)} = L^T \widetilde{\bm{n}}^{(k)}_i - 
        \widetilde{\bm{m}}_i$\;
        
        $W^{(k)}_i = \diag \left( \left( \left(\bm{u}^{(k)}\right)^2 + 
        \varepsilon^2 \right)^{\frac{\tau}{2} - 1} \right)$\;
        
        $\widetilde{\bm{n}}^{(k+1)}_i = \arg\min_{\widetilde{\bm{n}}} \left\| 
        \left(W_i^{(k)}\right)^{\frac{1}{2}} \left(L^T \widetilde{\bm{n}} - 
        \widetilde{\bm{m}}_i\right) \right\|_2^2$\; \label{mini}
        }
    
    \If{$\left\| \widetilde{N} - \widetilde{N}_{\rm old} \right\|_F < \gamma 
    \left\| \widetilde{N}_{\rm old} \right\|_F$}{
        \textbf{break}\;} \label{cond}
        }

\end{algorithm}

\subsection{Numerical tests}\label{subsec:refne}

Since our approach to deal with reflective surfaces is proposed here for the
first time, we report some numerical experiments to illustrate its performance.
We first analyze an abstract example, to test the method regardless or our
application. Then, we process a synthetic dataset
corresponding to a
realistic scene, constructed using the data generator described
in~\cite{cr2025}.

\paragraph{Ideal test}

We create an ideal test case composed of $q=10$ images, with $p=10000$ pixels, 
by setting random light directions in a matrix $L$ and a random normalized
vectors, in the matrix $N_{\text{true}}$, forming the normal field to a non
smooth surface.
The data matrix $M$ is computed as $(N_{\text{true}})^T L$, and a sparse
random error matrix $E$ is added to it. 
We let $\tau=0.1$ in the $\tau$-norm, $\varepsilon = 10^{-2}$ in~\eqref{Wk},
$\gamma=10^{-6}$ in~\eqref{stop}, and the maximum number of iteration is fixed
at $\textit{maxit}=100$.

To evaluate the solution, we compute the \emph{relative error}
\begin{equation}\label{RRE}
\text{RRE} (N)  = \frac{ \|N - N_{\text{true}} \|_F}{\| N_{\text{true}}\|_F},
\end{equation}
where $\|\cdot\|_F$ is the Frobenius norm and 
$N$ is the estimated normal vector field.
We also consider \emph{Type I} and \emph{Type II}  errors for the matrix
$E$,~\cite{buccini2023}, defined by

\begin{equation}\label{errtype}
\begin{aligned}
\text{Error Type I} (E) & = \frac{ \lvert \{ j: E_j \neq 0,
(E_{\text{exact}})_j = 0\} \rvert }{ \lvert \{ j: (E_{\text{exact}})_j = 0\} 
\rvert }, \\
\text{Error Type II} (E) & = \frac{ \lvert \{ j: E_j = 0, (E_{\text{exact}})_j
\neq 0\} \rvert }{ \lvert \{ j: (E_{\text{exact}})_j \neq 0\} \rvert},
\end{aligned}
\end{equation}
where $|\cS|$ denotes the number of elements of the set $\cS$.

The \emph{Type I} error corresponds the probability of incorrectly identifying
a zero element as a non-zero one. It is calculated as the ratio of zero entries
in $E$ that the algorithm incorrectly estimates as non-zeros to the total
number of true zero entries. In contrast, the \emph{Type II} error reflects the
probability of wrongly classifying a non-zero element as a zero. It is defined
as the ratio of non-zero entries in $E$ incorrectly sets to zero to the total
number of true non-zero entries. 

Row 1 of~\Cref{tab:errors} shows results obtained applying 
model~\eqref{reflprobl} to above data. 
We can observe that the method converges very fast. In just
$5$ iterations the RRE value drops to $10^{-7}$, so the
estimation of the normal vector field is very accurate. Moreover, Type I and
Type II errors are zero, indicating that all reflection points are correctly
detected by the method.

\begin{table}
\centering
\caption{Errors on the numerical tests. We consider the number of iterations
(\texttt{niter}), the RRE, and Type I and Type II errors.}
\label{tab:errors}
\begin{tabular}{l|l|l|l|l|l}
& \texttt{niter} & RRE & q & Type I & Type II\\ \hline
Ideal Test case & 5 & 9.21e-07 & all & 0 & 0 \\ \hline
\multirow{9}{*}{Synthetic dataset} &\multirow{9}{*}{28} & 
\multirow{9}{*}{3.97e-03} & all & 5.74e-04 & 2.55e-01 \\
& & & 1 & 7.04e-04 & 9.13e-02 \\
& & & 2 & 0 & 4.34e-01  \\
& & & 3 & 2.62e-03 & 2.50e-01 \\ 
& & & 4 & 0 & 8.43e-01 \\
& & & 5 & 0 & 9.13e-02 \\ 
& & & 6 & 0 & 8.19e-01 \\ 
& & & 7 & 1.31e-03 & 2.44e-01 \\ 
& & & 8 & 0 & 0
\end{tabular}
\end{table}

\paragraph{Synthetic dataset}

\begin{figure}
\centerline{
\includegraphics[width=.6\textwidth]{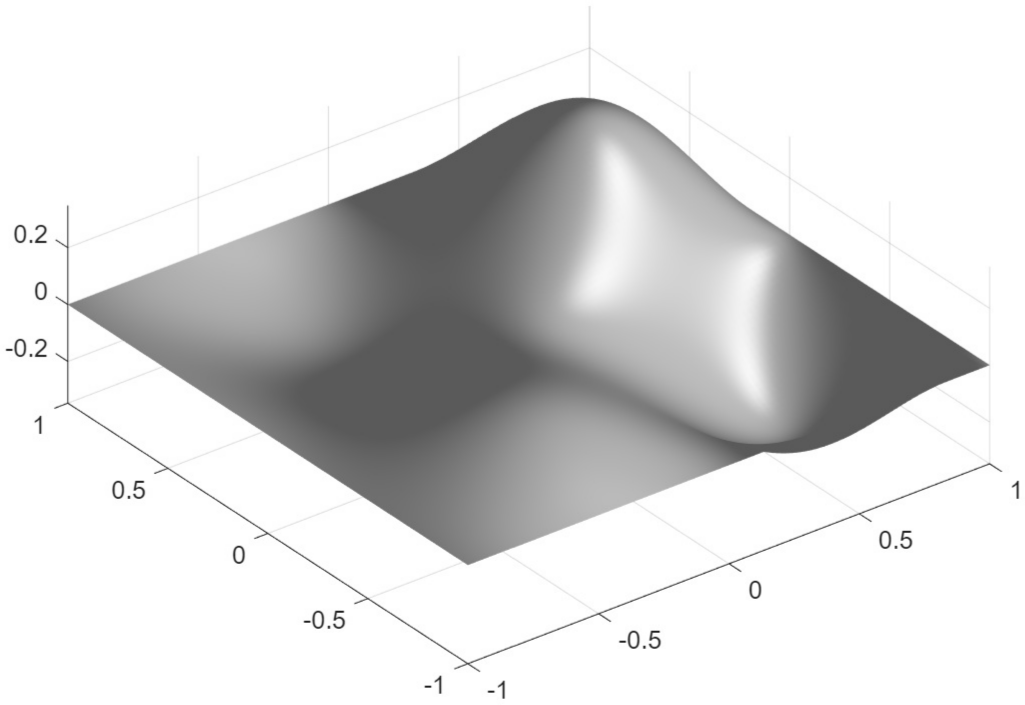}
\includegraphics[width=.4\textwidth]{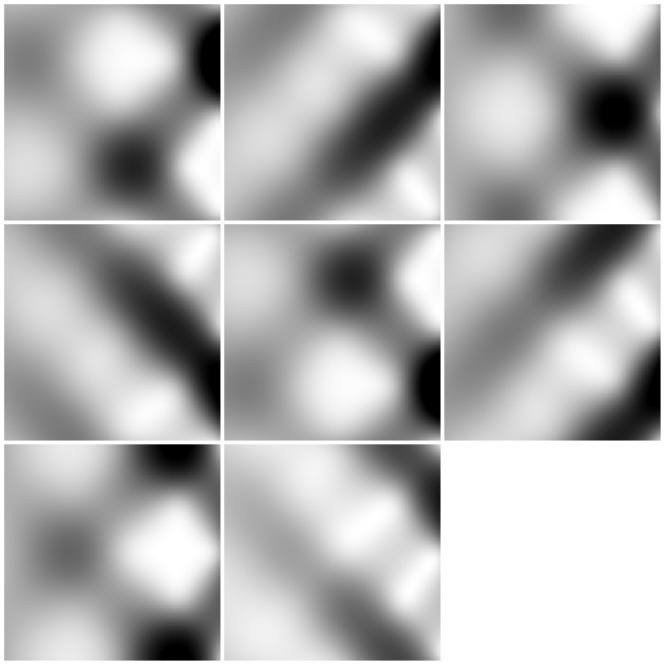} }
\caption{A synthetic surface used to study the performance of the model, and 
the corresponding dataset obtained considering $8$ different light directions.}
\label{fig:synthsurface}
\end{figure}

\begin{figure}
\centering
\subfigure[Dataset with reflection.]{
\includegraphics[width=.4\textwidth]{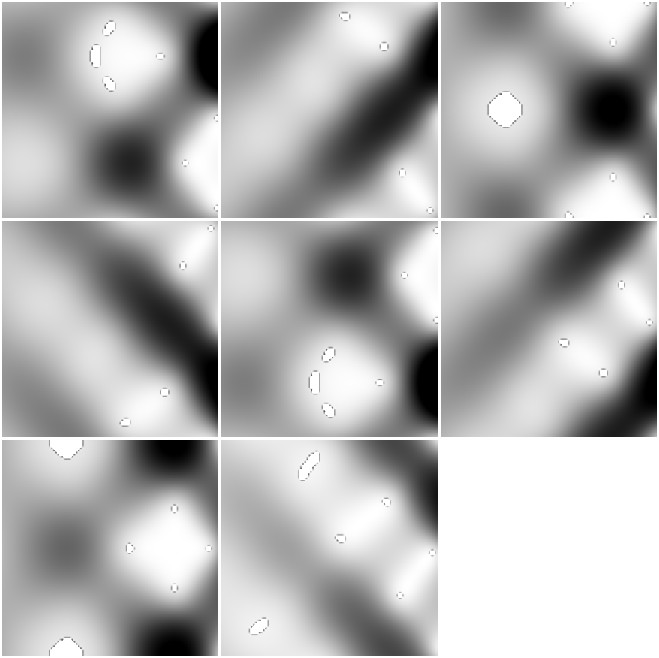} }
\subfigure[Reflection points.]{
\includegraphics[width=.4\textwidth]{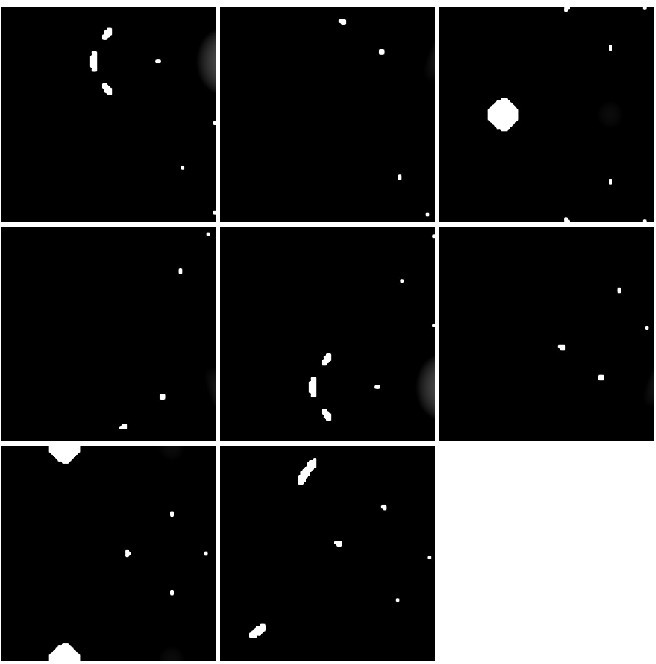}
\label{fig:reflsynth_pts} }
\caption{A synthetic dataset with reflection: images in the dataset (left) and
reflection points added to the data (right).}
\label{fig:reflsynth}
\end{figure}

Let us consider the smooth surface defined by
\begin{equation}\label{eqsynth}
u(x,y) = \frac{1}{5} e^x \sin{(\pi x)} \sin{(\pi y)}
\end{equation}
on the square domain $[-1, 1] \times [-1, 1]$, displayed
in~\Cref{fig:synthsurface} (left panel). 
\Cref{fig:synthsurface} (right panel) shows eight pictures of the surface with 
resolution $101 \times 101$ generated under fixed light positions,
specified in a matrix $L$.
The normal vectors, stored in the matrix $N$, are computed from partial
derivatives of the surface, and 
the matrix $E$ contains reflection points, simulating real reflection from the
surface in the direction of the observer.
This construction is performed by the package described
in~\cite{cr2025}. The final dataset and the reflection points are displayed
in~\Cref{fig:reflsynth}.
The parameters of the model are fixed as follows: $\tau =0.5$,
$\varepsilon=10^{-2}$, $\textit{maxit}=100$, and $\gamma =10^{-6}$.

\begin{figure}
\centering
\subfigure[Reconstruction.]{
\includegraphics[width=.6\textwidth]{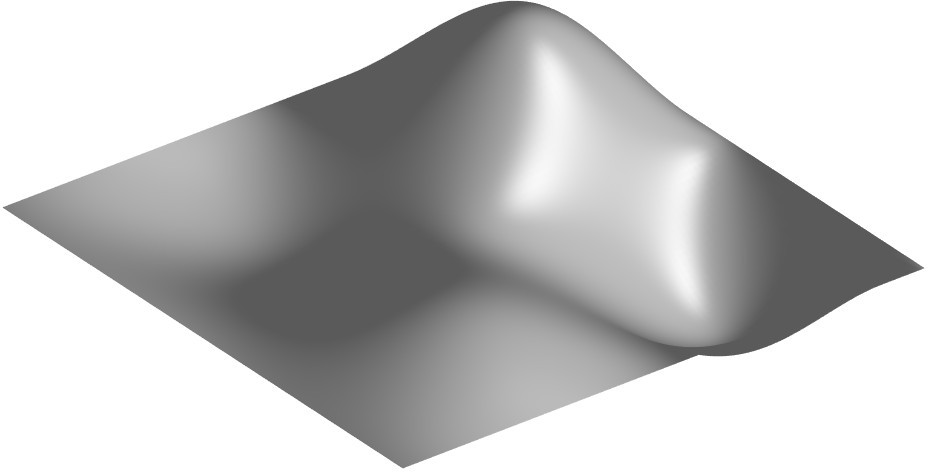}}
\subfigure[Detected reflection points.]{
\includegraphics[width=.35\textwidth]{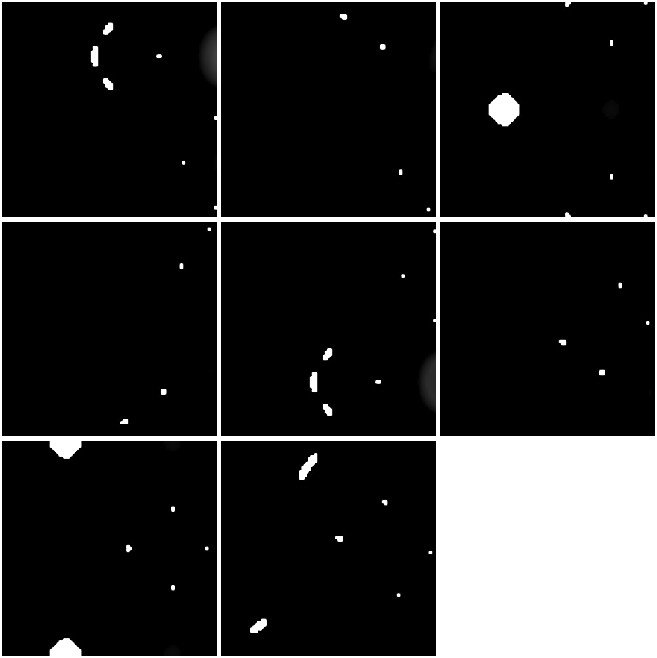}
\label{fig:reflsynth_detpts} }
\caption{Reconstructed surface and detected reflection points.}
\label{fig:reflsynthsol}
\end{figure}

Obtained results are illustrated in Row 2 of~\Cref{tab:errors}, where we
first report the errors across all the images, followed by Type I and Type
II errors for each image. In this case, after few iterations the method
estimates normal vectors with an error of $3.97 \cdot 10^{-3}$. 
Type I and Type II errors calculated across all images are consistently low,
with values around $10^{-4}$ and $10^{-1}$, respectively.
While Type II errors computed for individual images are generally in line
with the global value, in many cases Type I errors drop to zero.
The 3D reconstruction and the detected reflection points are reported in
\Cref{fig:reflsynthsol}.

We remark that solving this numerical problem by the standard PS method with
Hayakawa's estimation of the light directions results in a breakdown of the
algorithm. Neither one of the indicator discussed in~\Cref{sec:ideal} is
able to reduce the dataset to make the problem solvable.
If the matrix $L$ is assumed to be known and equation~\eqref{eqmatrix} is
solved directly, the relative error RRE for the normal matrix takes the value
3.82e-01, to be compared to the error 3.97e-03 reported in~\Cref{tab:errors}.

\section{The \pstredtools\ package}\label{sec:package}

An important aspect of this paper is the practical implementation of the models
discussed, which bridges theory and real-world applications. 

Here, we describe the \pstredtools\ package, built for this purpose. It
implements the photometric stereo model, enabling the reconstruction of the
shape of objects from multiple images, both grayscale or RGB.
Here we propose a revised version of the software first presented
in~\cite{cdfrv21}.
The toolbox is developed in the Matlab environment~\cite{matlab}, but it is
compatible with GNU Octave~\cite{octave}, an open source package with
Matlab-like syntax.
It is available at \url{https://bugs.unica.it/cana/software.html}.
The software has already been used for processing real archaeological
data; an example is documented in~\cite{cprt25}.

The \pstredtools\ package is composed of $20$ functions, $5$ datasets, and
$8$ demonstration scripts, listed in Table~\ref{table:ps3dfun}. 
Additionally, a
directory \emph{data} is provided to assist users in data preparation.
Lastly, the \emph{ps3dtools\_manual} PDF file contains a detailed
documentation for the package. Here we report a brief description.

\begin{longtable}{ l p{9.4cm}}
\caption{Functions and scripts in the \pstredtools\ package.}
\label{table:ps3dfun} \\
\hline
\multicolumn{2}{c}{\textbf{Reconstruction functions}} \\*
\ttt{ps3d} & photometric stereo 3D reconstruction \\
\ttt{ps3dstep1} & albedo, normals and lights reconstruction (first step) \\
\ttt{ps3dstep2} & returns the 3D coordinates of the object (second step) \\
\hline
\multicolumn{2}{c}{\textbf{Ideality indicators}} \\*
\ttt{psideal1} & extracts a subset of images based on ideality indicator \#1 \\
\ttt{psideal2} & extracts a subset of images based on ideality indicator \#2 \\
\hline
\multicolumn{2}{c}{\textbf{Reflection model}} \\*
\ttt{refl} & detects reflection points, estimate the normal vector field \\
\ttt{typeerrors} & computes Type I and Type II errors for the residual matrix E 
\\
\hline
\multicolumn{2}{c}{\textbf{Data generation}} \\*
\ttt{datagen} & synthetic dataset generator \\
\hline
\multicolumn{2}{c}{\textbf{Auxiliary routines}} \\*
\ttt{cholspec} & computes the Cholesky factor via SVD forcing positivity \\
\ttt{lightnlfun} & nonlinear function and Jacobian for ideality indicator \#2 \\
\ttt{normlights} & normalizes columns of L (Hayakawa's procedure) \\
\ttt{plotlights} & 2-D plot of light directions \\
\ttt{plotlights3d} & 3-D plot of light directions \\
\ttt{poisdir} & Poisson solver, Dirichlet boundary conditions \\
\ttt{poisneu} & Poisson solver, Neumann boundary conditions \\
\ttt{rotatelights} & finds the right rotation of the recovered lights \\
\ttt{setoptions} & sets default options for the \ttt{ps3d} function \\
\ttt{showdset} & displays a ps3d dataset \\
\ttt{tmngn} & Gauss-Newton method regularized by TSVD~\cite{pr25} \\
\ttt{writeply} & converts .mat surface into .ply file \\
\hline
\multicolumn{2}{c}{\textbf{Demonstration scripts and datasets}} \\*
\ttt{testshell} & demo test for the shell dataset \\
\ttt{testsynth} & demo test for a grayscale synthetic dataset \\
\ttt{testsynthcol} & demo test for a color synthetic dataset \\
\ttt{testsynthrefl} & demo test for a synthetic dataset with reflection \\
\ttt{testideal} & demo test for ideality indicators \\
\ttt{testdata} & demo test for a user supplied dataset \\
\ttt{showshell} & displays the shell dataset. \\
\ttt{showsynth} & displays the synthetic dataset \\
\ttt{shell.mat} & seashell dataset \\
\ttt{synth.mat} & synthetic dataset \\
\ttt{synthideal.mat} & synthetic dataset to test ideality indicators \\
\ttt{synthcol.mat} & colored synthetic dataset \\
\ttt{synthrefl.mat} & synthetic dataset with reflection points \\
\hline
\multicolumn{2}{c}{\textbf{Dataset construction support (subdirectory 
\ttt{data})}} \\*
\ttt{createdataset} & creates dataset for the \ttt{ps3d} function \\
\ttt{createdatasetrgb} & creates an RGB dataset for the \ttt{ps3d} function \\
\ttt{findoffset} & determines the alignment of the images \\
\ttt{png2ppm} & converts PNG file into PPM \\
\ttt{ppm2pgm} & converts PPM file into PGM \\
\hline
\end{longtable}

\subsection{Reconstruction functions}

This set of functions implements the photometric stereo model under the 
assumption of unknown lighting. The main function, \ttt{ps3d}, relies on two 
key steps: 
\begin{itemize}
\item \ttt{ps3dstep1} computes the normal vector field by leveraging
\ttt{normligths}, that estimates the normal vector field $N$ and the light
directions $L$, starting from the initial estimation obtained from the SVD of
$M$, and \ttt{rotatelights}, that returns the final estimates of the matrices
$N$ and $L$, and the rotation matrix $Q$;
\item \ttt{ps3dstep2} determines the 3D coordinates of the observed object,
depending on Dirichlet or Neumann boundary conditions. 
\end{itemize}

Some parts of the computation are included in a set of \emph{Auxiliary
functions}; see Table~\ref{table:ps3dfun}.
Such group also includes various data visualization
routines, such as \ttt{plotlights} and \ttt{plotlights3d}, to represent 
light directions, \ttt{showdset}, to display a dataset, and \ttt{writeply},
which stores a computed surface in a \ttt{.ply} file. Such format is
suitable to drive 3D printers and to visualize the reconstructed surface on a
3D viewer/editor like, e.g., MeshLab~\cite{meshlab}.
Unfortunately, \ttt{writeply} cannot be used in Octave.

\subsection{Ideality indicators}

This pair of functions implements the computation of the ideality
indicators described in~\cite{cprt23}. 
They extract from a dataset a subset of images that
meet the specific condition considered by the corresponding indicator.
The first indicator of ideality (\ttt{psideal1}) tries to ensure
the positiveness of the matrix $G$ in~\eqref{hayaprob}, while the second one
(\ttt{psideal2}) aims at keeping the Jacobian of~\eqref{nonlineq} full rank;
see \Cref{sec:ideal}.

\subsection{Reflection model}

The function \ttt{refl} implements the model introduced
in~\Cref{sec:reflmod} to process data affected by reflection
phenomena. It detects reflection points in data and estimates the normal
vector field. 
Type I and Type II errors for the reflection matrix $E$ are computed by
the \ttt{typeerrors} function.

\subsection{Data generation}

A data generation script is useful for the user to make experiments
under specific chosen conditions. 
The function \ttt{datagen}, by fixing certain parameters, creates a dataset
that reflects desired observation conditions. 
Such function is strictly related to the
\ttt{psdatasynth} package, described in~\cite{cr2025}.
It is, in fact, a simplified version of the \ttt{psdatagen} script from that
package.
If \ttt{psdatagen} is installed, \ttt{datagen} invokes some of its functions 
to generate a set of synthetic pictures. Otherwise, some predefined datasets are
loaded.
Demo tests available in \pstredtools\ do not require the presence of the
data generation package.

\subsection{Demonstration scripts and datasets}

Four synthetic datasets, generated using the \ttt{psdatasynth}
toolbox~\cite{cr2025}, are included in this package. 
There is also an experimental dataset that contains pictures of a shell under
direct sun light, so that light rays can be considered virtually parallel; 
see~\cite{cdfrv21} and~\cite{rrmyb23} for discussions on this dataset.

Each dataset comes with a demonstration script that shows how to apply the PS
model to the corresponding set of images.
Some tests also allow to modify the observation conditions, if the
\ttt{psdatagen} toolbox is available.
Additionally, two visualization scripts,
\ttt{showshell} and \ttt{showsynth}, are included to demonstrate how a dataset
can be visualized.

\paragraph{Datasets}
\begin{itemize}

\item {\ttt{shell}}: experimental dataset composed of $20$ images of a real
shell under sunlight, displayed in \Cref{fig:shell}. 

\begin{figure}
\centerline{
\includegraphics[width=.6\textwidth]{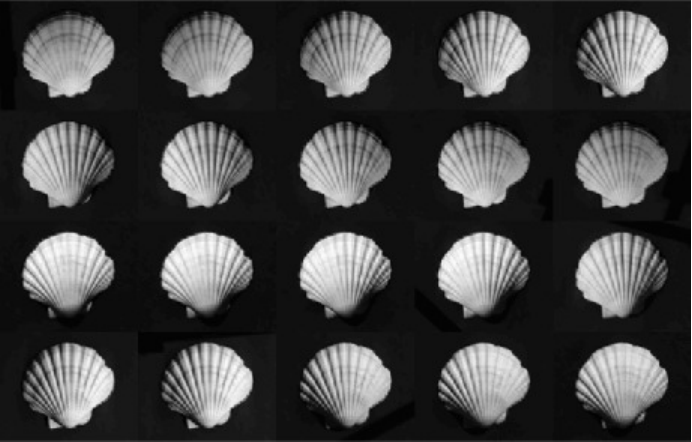} }
\caption{The $20$ images in the \ttt{shell} experimental dataset; 
see~\cite{cdfrv21,rrmyb23} on its creation.}
\label{fig:shell}
\end{figure}

\item {\ttt{synth}}: basic synthetic dataset, shown
in~\Cref{fig:synthsurface}.

\item {\ttt{synthcol}}: color (RGB) synthetic dataset, obtained from the
{\ttt{synth}} dataset, using the \ttt{mandrill} demo image from Matlab as the
albedo matrix; it is displayed in~\Cref{fig:synthmandrill}.

\begin{figure}
\centering
\includegraphics[width=.6\textwidth]{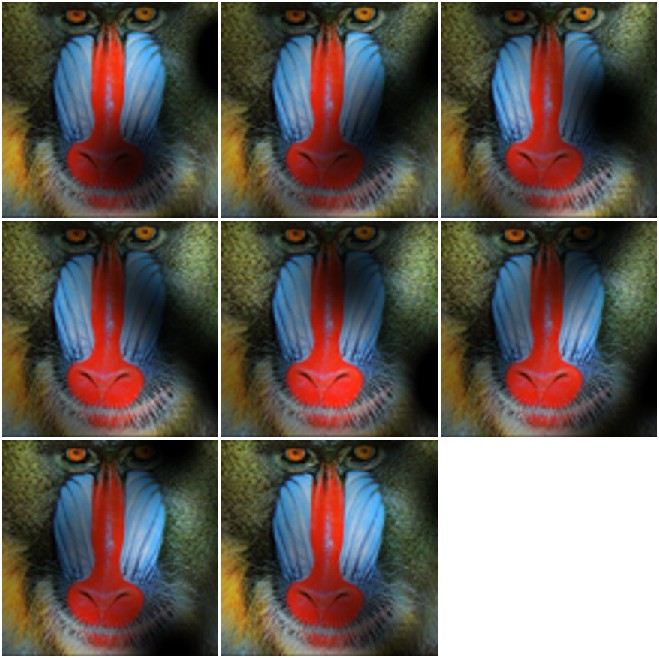}
\caption{Colored synthetic dataset obtained using surface~\eqref{eqsynth} with 
the \ttt{mandrill} image as albedo.}
\label{fig:synthmandrill}
\end{figure}

\item {\ttt{synthrefl}}: it consists of the {\ttt{synth}} dataset with 
reflection points added consistently with the shape of the surface;
see~\Cref{fig:reflsynth}.

\item {\ttt{synthideal}}: synthetic dataset to test ideality indicators.
Two data matrices are included: \ttt{Minf} contains 9 images obtained
using different light sources placed at infinite distance from a surface;
\ttt{MM3} is composed of images obtained from light sources at distance 
$\delta\cdot A$, where $A$ is the size of the observed object and 
$\delta=10,8,6,4,2$; see the \ttt{testideal} script below.

\end{itemize}

\paragraph{Demo tests}

\begin{itemize}
\item {\ttt{testshell}}:
demo test for the shell dataset. The results obtained applying the \ttt{ps3d}
function to the \ttt{shell} dataset are displayed in~\Cref{fig:shellrec}.

\begin{figure}
\centerline{
\includegraphics[width=.6\textwidth]{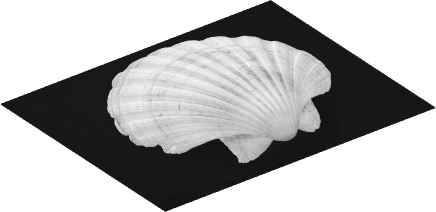} }
\caption{Reconstructed shell surface.}
\label{fig:shellrec}
\end{figure}

\item {\ttt{testsynth}}:
demo test for a grayscale synthetic dataset. 
\Cref{fig:synthrec} shows the reconstruction obtained processing the
\ttt{synth} dataset.

\begin{figure}
\centerline{
\includegraphics[width=.6\textwidth]{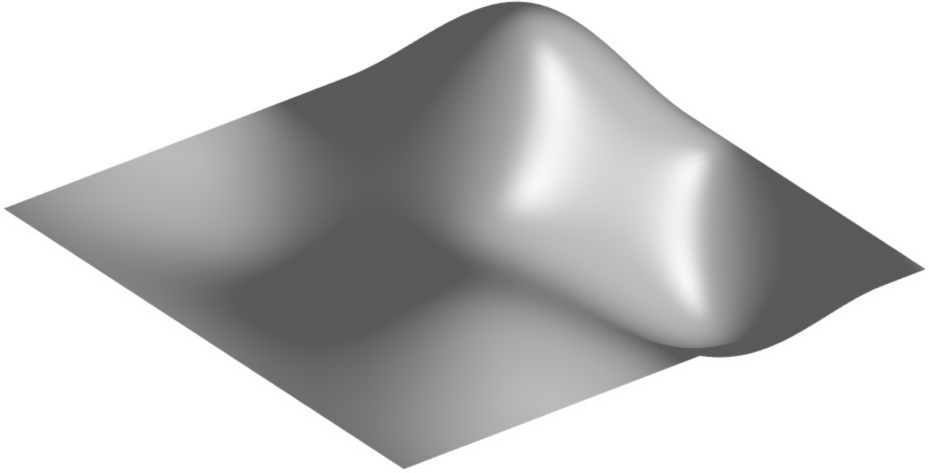} }
\caption{The restored surface obtained applying the photometric stereo
technique to the \ttt{synth} dataset.}
\label{fig:synthrec}
\end{figure}

\item {\ttt{testsynthcol}}:
demo test for a color synthetic dataset. The data file
\ttt{synthcol} is processed. The blue channel is used to determine the
shape of the object, while the other two channels only to estimate the remaining
color information. The reconstruction is displayed in~\Cref{fig:synthcolrec}. 

\begin{figure}
\centerline{
\includegraphics[width=.6\textwidth]{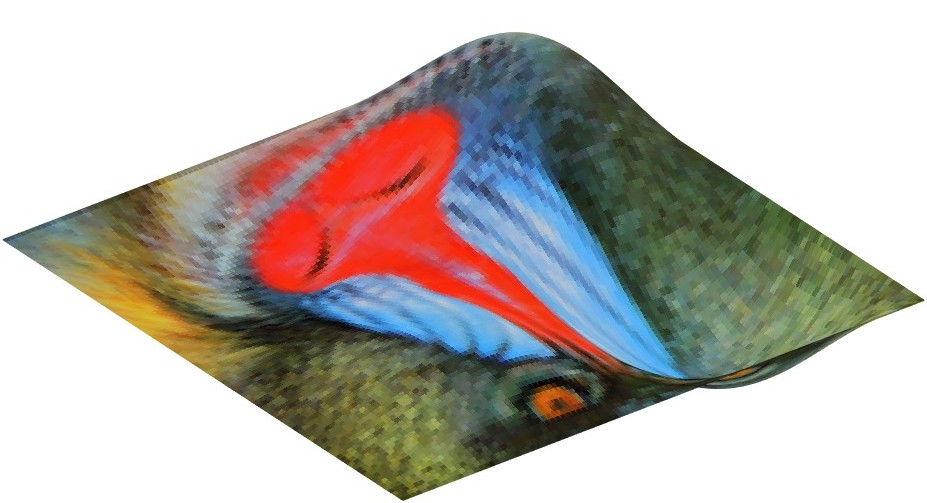} }
\caption{Color reconstruction of the \ttt{synthcol} surface obtained using the 
PS model.}
\label{fig:synthcolrec}
\end{figure}

\item {\ttt{testsynthrefl}}:
demo test for a synthetic dataset that presents reflections. After
loading the \ttt{synthrefl} dataset, the \ttt{refl} function is employed to
estimate the normal vector field $N$. The relative error~\eqref{RRE} on the 
normal vectors, together with Type I and Type II errors~\eqref{errtype}, are
computed. Finally, the function \ttt{ps3dstep2} is used to determine a
reconstruction of the surface starting from the estimated normal vector field;
see Figures~\ref{fig:reflsynth} and~\ref{fig:reflsynthsol}.

\item {\ttt{testideal}}:
demo test for ideality indicators.
Once the \ttt{synthideal} dataset is loaded, the third image is
replaced with a picture of the same object lighted by a source placed at
different close positions. 
This process is used to evaluate the capability of ideality indicators to
recognize the distorted image. The
relative errors on the reconstruction corresponding to the original dataset
and the reduced ones are computed and displayed in~\Cref{fig:indres}.

\begin{figure}
\centerline{
\includegraphics[width=.6\linewidth]{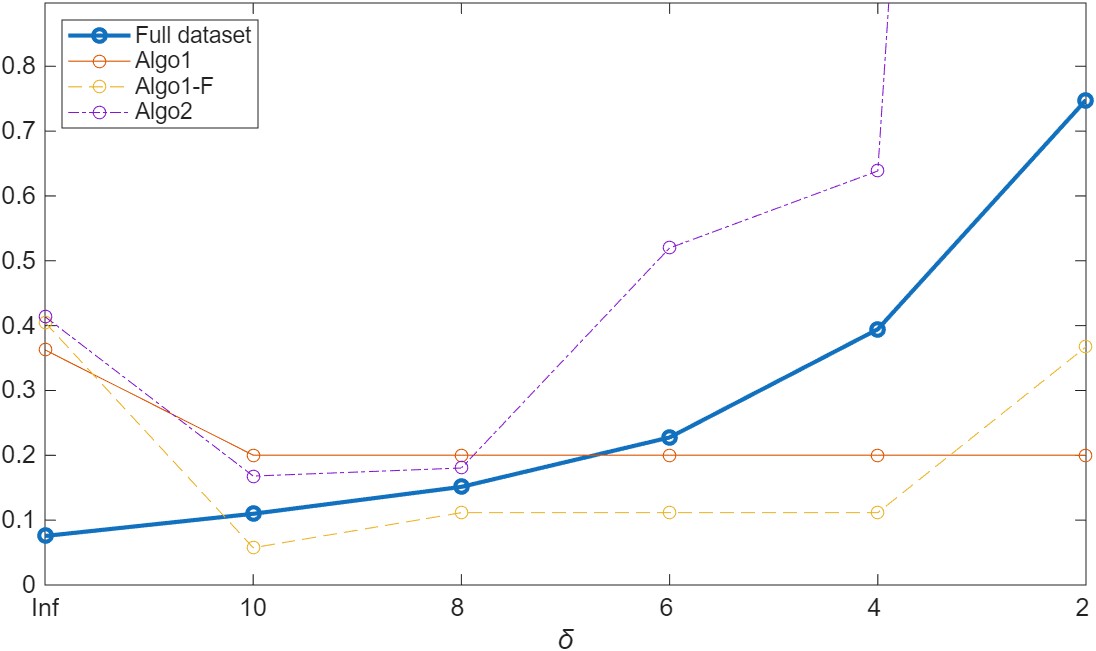}}
\caption{Relative errors in the $\infty$-norm, obtained by processing either 
the whole dataset or the reduced datasets, obtained applying both indicators of 
ideality; see \Cref{sec:ideal}. $\delta A$ 
is the distance of the light source from the object.}
\label{fig:indres}
\end{figure}

\item {\ttt{testdata}}:
demo test for a user supplied dataset. The dataset has to be prepared
according to the format used in the package. The functions in the 
\emph{data} directory may help the user in doing that.

\end{itemize}

\subsection{Dataset construction support}

The \ttt{data} directory contains a set of functions that supports the user in
creating a dataset in the format required by the package, from a set of
acquired images.
The \ttt{findoffset} script is used to detect the offset for aligning a
dataset. Such process is motivated by small displacements characterizing
experimentally acquired images, and often enhances the definition of the
reconstructed surface.
The functions \ttt{createdataset} and \ttt{createdatasetrgb} can be
used to create a dataset from a collection of pictures.
Two conversion routines may be helpful to speed up data processing.

We mention the fact that digital cameras usually produce compressed images,
which may introduce distortion in 3D reconstructions.
Modern cameras allow to store images in ``RAW'' format, avoiding any
processing by the camera software.
Such ``raw'' images are memory consuming, since they are not compressed, and
are often stored using a proprietary format which is not supported by common
software packages.
To convert them to a suitable format (e.g., PNG) a specific computer program
is required. We advise the use of the open source application Raw
Therapee~\cite{rawtherapee} for this task.

\section{Conclusion}\label{sec:con}

In this paper, we presented a revised and extended version of the Matlab package
\ttt{\small ps3d}, designed to support the photometric stereo model for the 3D
reconstruction of real objects from images acquired under varying illumination
conditions.

The new release addresses several limitations of the original package. In
particular, we introduce an alternative model to handle reflective surface
phenomena, which represent one of the most common sources of distortion in
real-world acquisitions. Furthermore, the package now includes dedicated
indicators to quantitatively assess the compliance of a given dataset with the 
PS model assumptions.
Compatibility with GNU Octave has also been ensured.
The tool supports both synthetic and experimental data, and it is completed by 
some demo tests.

\section*{Acknowledgment}

The authors are members of the GNCS group of INdAM.
A.B. and F.P. are partially supported by the INdAM-GNCS project 2026 ``Metodi 
strutturati per il signal processing avanzato'' (CUP E53C25002010001). 
G.R. is partially supported by the INdAM-GNCS project 2026 ``Metodi 
numerici per modelli integrali e dinamiche con memoria'' (CUP E53C25002010001).
A.B. and G.R. are partially supported by the research project Fondazione di
Sardegna 2024-2025 ``Integral and Discrete Inverse Problems (InDIP)'' (CUP
F83C26000350007).
F.P. gratefully acknowledges the financial support by Regione Lombardia under
the PR FESR 2021-2027 - Collabora \& Innova Project ``ID-Tex'' (CUP
E19I25000950007).

\printbibliography

\end{document}